\documentclass[12pt]{article}
\title{\bf Some characterizations of freeness of hyperplane arrangements}
\author{
\textsc{Masahiko Yoshinaga} \\[.1in]
\texttt{yosinaga@kurims.kyoto-u.ac.jp}\\[.1in]
Research Institute for Mathematical Sciences\\
Kyoto University\\
Kyoto, 606-8502 Japan
}
\date{\today}

\usepackage{amssymb}
\usepackage{amscd}
\usepackage[all]{xy}

\newtheorem{Def}{Definition}[section]
\newtheorem{Prop}[Def]{Propositon}
\newtheorem{Thm}[Def]{Theorem}
\newtheorem{Lemma}[Def]{Lemma}
\newtheorem{Cor}[Def]{Corollary}
\newtheorem{Rem}[Def]{Remark}

\newtheorem{Example}[Def]{Example}

\newtheorem{Assum}[Def]{Condition}
\newtheorem{Conjecture}[Def]{Conjecture}
\newtheorem{Conclusion}[Def]{Conclusion}

\newcommand{\bbK}{\mathbb{K}}
\newcommand{\bbP}{\mathbb{P}}
\newcommand{\bbR}{\mathbb{R}}

\newcommand{\bbZ}{\mathbb{Z}}

\newcommand{\bfk}{\mathbf{k}}
\newcommand{\bfres}{\mathbf{res}}

\newcommand{\rmH}{\mathrm{H}}

\newcommand{\calA}{\mathcal{A}}

\newcommand{\calE}{\mathcal{E}}
\newcommand{\calF}{\mathcal{F}}

\newcommand{\calO}{\mathcal{O}}

\newcommand{\cone}{\mathbf{c}}
\newcommand{\rmExt}{\mathrm{Ext}}
\newcommand{\calExt}{\mathcal{E}xt}

\newcommand{\rmDer}{\mathrm{Der}}

\newcommand{\del}{\partial}
\newcommand{\til}{\widetilde}
\newcommand{\owari}{\hfill$\square$\medskip}

\begin{document}

\maketitle


\begin{abstract}
We consider some characterizations of freeness of 
a hyperplane arrangement, in terms of the following 
properties: 
local freeness, factorization of the characteristic polynomial, and 
freeness of the restricted multiarrangement. 
In the case of a $3$-arrangement, freeness is characterized by 
factorization of the characteristic polynomial and coincidence 
of its roots 
with the exponents of the restricted multiarrangement. 
In the case of higher dimensions, 
it is characterized by a kind of local freeness and 
freeness of the restricted multiarrangement. 
As an application, we prove the freeness of certain 
arrangements conjectured by Edelman and Reiner. 
\end{abstract}

\section{Introduction}
\label{sec:intro}
A hyperplane arrangement $\calA$ in $\ell$ dimensional linear space 
is said to be free with exponents $\exp(\calA)=(1=d_1, d_2, \cdots, d_\ell)$ if
the associated module of all logarithmic vector fields is free 
with basis $\delta_1, \cdots, \delta_\ell$ such that 
$\deg \delta_i =d_i$. 
Among other properties we consider the following three necessary 
conditions for an arrangement to be free, which are due to Terao and Ziegler: 
\begin{itemize}
\item[(1)] $\calA$ is locally free. 
\item[(2)] The characteristic polynomial $\chi(\calA, t)$ factors 
completely over $\bbZ$, indeed, it is equal to 
$(t-d_1)(t-d_2)\cdots (t-d_\ell)$. 
\item[(3)] The multiarrangement obtained by restricting to 
a hyperplane is free with multiexponent 
$(d_2, \cdots, d_\ell)=\exp(\calA)\backslash\{1\}$. 
\end{itemize}
It is known that each condition is not sufficient to characterize 
freeness. For example, since any central $2$-(multi)arrangement is free, 
$3$-arrangements are always locally free and restricted multiarrangements 
are also free. So (1) and (3) hold 
for any $3$-arrangement. But it is not necessarily free. 
We also note that Kung \cite{kun} found many examples of non-free 
arrangements whose characteristic polynomial factors completely over 
$\bbZ$. Recently, Schenck \cite{sch-recent} studies 
the difference between freeness and the factorization of characteristic 
polynomials for $3$-arrangements. 

It seems natural to ask whether a combination of some conditions 
characterizes freeness. 
The behavior is completely different for $\ell =3$ or $\ell \geq 4$. 

For $3$-arrangement $\calA$ 
we will prove, in \S\ref{sec:3arr}, that $\calA$ is free 
if and only if it satisfies (2)$+$(3) with the coincidence of numbers, 
i.e. the condition 
(2+3) below characterizes freeness. (Theorem \ref{thm:3})
\begin{itemize}
\item[(2+3)] The characteristic polynomial factors as 
$$
\chi(\calA, t)=(t-1)(t-d_2)(t-d_3) 
$$
and multiexponents of restricted multiarrangement is 
$(d_2, d_3)$. 
\end{itemize}
Our proof is based on a study of Solomon-Terao's formula 
for characteristic polynomial 
and Ziegler's restriction map using Hilbert series. 

In \S\ref{sec:4}, we will study the freeness of an arrangement 
from the viewpoint of coherent sheaves on projective space. 
The freeness of an $\ell(\geq 4)$-arrangement can be 
characterized by (1) and (3). Furthermore, with the help of 
results on reflexive sheaves, 
we will give a characterization 
by the following weaker condition (Theorem \ref{thm:FREE}). 
\begin{itemize}
\item[(1'+3)] Let $H\in\calA$ be a hyperplane. 
The restricted multiarrangement 
on $H$ is free and $\calA$ is locally free along $H$. 
\end{itemize}

Edelman-Reiner \cite{ede-rei} conjectured that cones over certain 
truncated affine Weyl arrangements are free. 
As an application of our characterization of freeness, 
we prove that the Edelman-Reiner conjecture is true. 
The first half of condition (1'+3) has been proved by Terao \cite{ter-multi}. 
The second half will be proved by induction on the rank of root system 
using following fact: 
any localization of Weyl arrangement decomposes into a 
direct sum of Weyl arrangements of lower ranks. 
Indeed, the required local freeness is equivalent to 
the Edelman-Reiner conjecture for 
root systems of strictly lower ranks. The problem is resolved into 
computation of characteristic polynomials for 
rank two root systems, which has been verified by Athanasiadis 
\cite{ath-adv, ath-def}. 
And we also give a family of free arrangements which interpolates 
between extended Shi and extended Catalan arrangements. 

\noindent
{\bf Acknowledgement.} 
The author is grateful to Professor H. Terao for pointing out 
a mistake in the first draft of this paper and giving 
a correct argument. 
The author is grateful to Professor K. Saito for his advice 
and support. The author also thanks many friends who kindly 
guided the author to basic facts on algebraic geometry and vector bundles. 

\section{Preliminaries}

Let $V$ be an $\ell$-dimensional linear space over an arbitrary field 
$\bbK$ of characteristic $0$ and 
$S:=\bbK[V^*]$ be the algebra of polynomial functions on $V$ that is 
naturally isomorphic to $\bbK[z_1, z_2, \cdots, z_\ell]$ for 
any choice of basis $(z_1, \cdots, z_\ell)$ of $V^*$. 

A (central) hyperplane arrangement $\calA$ is a finite 
collection of codimension one 
linear subspaces in $V$. 
For each hyperplane $H$ of $\calA$, 
fix 
a nonzero linear form $\alpha_H\in V^*$ vanishing on $H$ and 
put $Q:=\prod_{H\in\calA}\alpha_H$. 

The characteristic polynomial of $\calA$ is defined as 
$$
\chi(\calA, t)=\sum_{X\in L_\calA}\mu(X)t^{\dim X}, 
$$
where $L_\calA$ is a lattice consists of the intersections of 
elements of $\calA$, ordered by reverse inclusion, 
$\hat{0}:=V$ is the unique minimal element of $L_\calA$ and 
$\mu:L_\calA\longrightarrow\bbZ$ is the M\"obius function
defined as follows:
\begin{eqnarray*}
\mu(\hat{0})&=&1,\\
\mu(X)&=&-\sum_{Y<X}\mu(Y),\ \mbox{if}\ \hat{0}<X. 
\end{eqnarray*}

Denote by 
$\rmDer_V$ and $\Omega_V^p$, respectively, the $S$-module of 
all polynomial vector fields and of polynomial differential 
$p$-forms over $V$. 

\begin{Def}
\label{def:log}
For a given arrangement $\calA$, 
we define modules of logarithmic vector fields and 
logarithmic $p$-forms by, respectively, 
$$
D(\calA)=\{ \delta\in\rmDer_V\ |\ \delta(\alpha_H)\in \alpha_H S,\ 
\forall H\in\calA\}
$$
and 
$$
\Omega_V^p(\calA)=\left\{\left. \omega\in\frac{1}{Q}\Omega^p_V\ \right|\ 
Q\cdot\frac{d\alpha_H}{\alpha_H}\wedge\omega\in\Omega_V^{p+1},
\ \forall H\in\calA \right\}. 
$$
\end{Def}

Next we recall a formula due to Solomon and Terao 
which deduces 
the characteristic polynomial from the Hilbert series of 
the graded $S$-modules $\Omega^p(\calA)$. 
For a finitely generated graded $S$-module $M$, the series 
$P(M,x)\in\bbZ[x^{-1}][[x]]$ defined by 
$$
P(M,x)=\sum_{p\in\bbZ}(\dim_\bbK M_p)x^p 
$$
is called the Hilbert series.

For an arrangement $\calA$, define 
$$
\Phi(\calA ; x, y)=
\sum_{p=0}^\ell P(\Omega^p(\calA), x)y^p. 
$$

\begin{Thm}
{\normalfont \cite{sol-ter}}
\label{thm:st}
The characteristic polynomial $\chi(\calA, t)$ is expressed as 
$$
\chi(\calA, t)=\lim_{x\rightarrow 1}\Phi(\calA; x, t(1-x)-1). 
$$
\end{Thm}

An arrangement $\calA$ is said to be {\bf free} if 
$D(\calA)$ (or equivalently $\Omega^1(\calA)$) 
is a free $S$-module, and then the multiset of 
degrees $\exp(\calA):=(d_1, d_2, \cdots, d_\ell)$ 
of a homogeneous basis of $D(\calA)$ 
is called the {\bf exponents}. 

Theorem \ref{thm:st} 
yields a famous factorization theorem by Terao; 
\begin{Thm}
{\normalfont \cite{ter-free}}
\label{thm:fact}
If $\calA$ is a free arrangement with exponents $(d_1, d_2, \cdots, d_\ell)$, 
the characteristic polynomial factors as follows; 
$$
\chi(\calA, t)=\prod_{i=1}^\ell (t-d_i). 
$$
\end{Thm}

A multiarrangement (introduced by Ziegler \cite{zie}) 
is a pair $(\calA, \bfk)$ 
consisting of an ordinary arrangement  $\calA$ and a map 
$\bfk:\calA\rightarrow\bbZ_{\geq 0}$ 
(called multiplicity). 
Any arrangement can be considered as a multiarrangement 
with constant multiplicity $\bfk(H)=1$, $\forall H\in\calA$. 

\begin{Def}
\label{def:multilog}
Let $(\calA, \bfk:\calA\rightarrow\bbZ_{\geq 0})$ be a multiarrangement. 
Denote $Q(\calA, \bfk):=\prod_{H\in\calA}\alpha_H^{\bfk(H)}$. 
We define the modules $D(\calA, \bfk)$ and $\Omega^p(\calA, \bfk)$ by 
$$
D(\calA, \bfk)=\{ \delta\in\rmDer_V\ |\ 
\delta(\alpha_H)\in \alpha_H^{\bfk(H)} S,\ 
\forall H\in\calA\}
$$
and
$$
\Omega^p(\calA, \bfk)=\left\{\left. \omega\in\frac{1}{Q(\calA, \bfk)}\Omega^p_V\ 
\right|\ 
Q(\calA, \bfk)\cdot\frac{d\alpha_H}{\alpha_H^{\bfk(H)}}\wedge\omega\in\Omega_V^{p+1},
\ \forall H\in\calA \right\}. 
$$
\end{Def}
The following is straightforward. 
$$
\Omega^0(\calA, \bfk)=S,\ \ \Omega^\ell(\calA, \bfk)=\frac{1}{Q(\calA, \bfk)}
\Omega_V^\ell. 
$$
A multiarrangement $(\calA, \bfk)$ is said to be free if 
$D(\calA, \bfk)$ (or equivalently $\Omega^1(\calA, \bfk)$) 
is a free $S$-module. In this case, the multiset of 
degrees of a homogeneous basis of $D(\calA, \bfk)$ 
is called the {\bf multiexponents} 
and also denoted by $\exp(\calA, \bfk)$. We list some 
consequences of the freeness of a multiarrangements 
that will be used later. 

\begin{Thm}
\label{thm:saito's}
{\normalfont \cite{sai-log}(Saito's criterion)}
Let $\omega_1,\cdots, \omega_\ell\in\Omega^1(\calA, \bfk)$ be 
homogeneous and linearly independent over $S$. Then $(\calA, \bfk)$ is 
free with basis $\omega_1, \cdots, \omega_\ell$ if and only if 
$$
\sum_{i=1}^{\ell}\deg \omega_i =-\sum_{H\in\calA}\bfk(H). 
$$
\end{Thm}

\begin{Thm}
\label{thm:localization}
{\normalfont (Localization Theorem)}
Let $(\calA, \bfk)$ be a free multiarrangement. For any 
intersection $X\in L_\calA$, define 
$$
\calA_X:=\{ H\in\calA\ |\ H\ni X\ \}. 
$$
Then the induced multiarrangement $(\calA_X, \bfk|_{\calA_X})$ 
is also free. 
\end{Thm}

\begin{Example}
\label{ex:naturalmulti}
Let $\calA$ be an arrangement in $V$ and  $H_1\in\calA$ be a hyperplane. 
Then the restriction of $\calA$ to $H_1$ is the arrangement 
$\calA^{H_1}:=\left\{H\cap H_1\ |\ H\in\calA \backslash \{H_1\}\ \right\}$ 
in $H_1$. 
This restriction has a natural structure of a multiarrangement 
$(\calA^{H_1}, \bfk_\calA^{H_1})$ with multiplicity 
$\bfk_\calA^{H_1}:\calA^{H_1}\rightarrow\bbZ_{\geq0}$ 
defined by 
$$
\bfk_\calA^{H_1}:\calA^{H_1}\ni H'\longmapsto 
\sharp\left(\{ H\in\calA\ |\ H\cap H_1=H'\}\right). 
$$
\end{Example}
Fix a coordinate system  $(z_1, z_2, \cdots, z_\ell)$ 
such that  $H_1=\{z_1=0\}$. Then 
every $\omega\in\Omega^p(\calA)$ can be uniquely expressed as 
$$
\omega=\omega_1+\frac{dz_1}{z_1}\wedge\omega_2, 
$$
where $\omega_1$ and $\omega_2$ are rational differential 
forms in $dz_2, \cdots, dz_\ell$. 
The following theorem was proved by Ziegler \cite{zie}. 
(For another formulation by using the modules of vector fields, 
see Theorem \ref{thm:zie2}.) 
\begin{Thm}
\label{thm:zie}
Using the notation above, $\omega_1|_{H_1}$ is contained in 
$\Omega^p(\calA^{H_1}, \bfk)$. 
Furthermore, $\calA$ is free if and only if the restricted multiarrangement 
$(\calA^{H_1}, \bfk)$ is free and the restriction map (for $p=1$) 
$$
\begin{array}{ccc}
\Omega^1(\calA)& \longrightarrow & \Omega^1(\calA^{H_1}, \bfk)\\
&&\\
\omega=\omega_1+\frac{dz_1}{z_1}\wedge\omega_2&\longmapsto&\omega_1|_{H_1}
\end{array}
$$
is surjective. 
\end{Thm}
Let us denote by 
$$
\bfres_H^p:\Omega^p(\calA)\longrightarrow \Omega^p(\calA^H, \bfk)\ 
\left( \omega_1+\frac{dz_1}{z_1}\wedge\omega_2 \mapsto \omega_1|_H \right)
$$
the restriction map and by 
$M^p\subset\Omega^p(\calA^H, \bfk)$ 
the image of $\bfres_H^p$. 
$M^p$ is a graded $\bbK[H]:=S/z_1S$ submodule by definition. Though 
the next corollary is easily deduced from Theorem \ref{thm:zie} 
and Saito's criterion, it is 
a starting point of our characterization of freeness. 
\begin{Cor}
\label{cor:criterion}
If the restriction map 
$\bfres_H^1:\Omega^1(\calA)\rightarrow\Omega^1(\calA^H, \bfk)$ 
is surjective and the restricted multiarrangement 
$(\calA^H, \bfk)$ is free with multiexponents 
$(d_2', \cdots, d_\ell')$, then $\calA$ is free with 
exponents $(1, d_2', \cdots, d_\ell')$. 
\end{Cor}
{\bf Proof.} Surjectivity implies that there exist 
$\omega_2, \cdots, \omega_\ell\in\Omega^1(\calA)$ such that 
$\bfres_H^1(\omega_2),\cdots, \bfres_H^1(\omega_\ell)\in
\Omega^p(\calA^H, \bfk)$ form a basis. Then from Saito's criterion 
\ref{thm:saito's}, 
$\omega_2, \cdots, \omega_\ell$ together with 
$d\alpha/\alpha\in\Omega^1(\calA)$ form a basis of $\Omega^1(\calA)$.
\owari

Since $\Omega^\bullet(\calA)$ is closed under exterior product, 
we have a homomorphism (where $\alpha=\alpha_H$ is a defining equation 
of $H\in\calA$) 
$$
\del:\Omega^p(\calA)\longrightarrow\Omega^{p+1}(\calA),\ 
\left( \omega\longmapsto (d\alpha/\alpha)\wedge\omega\right). 
$$
\begin{Prop}
{\normalfont (\cite[Prop. 4.86]{orl-ter}) }
\label{prop:exact}
The complex $(\Omega^\bullet(\calA), \del)$ is acyclic. 
\end{Prop}
Using this proposition, we study the effect of 
restriction map on the Hilbert series. 
\begin{Thm}
\label{thm:hilbformula}
The Hilbert series of $M^p$ and $\Omega^p(\calA)$ are connected 
by the relationship 
\begin{equation}\label{eq:hilbformula}
\sum_{p=0}^{\ell-1}P\left(M^p, x\right)
y^p=
\frac{x(1-x)}{x+y}
\times\Phi(\calA; x,y).
\end{equation}
\end{Thm}

\noindent{\bf Proof.}
The restriction map above can be considered as 
a composition of $\frac{dz_1}{z_1}$ and the Poincar\'e 
residue map, 
$$
\xymatrix{
   \Omega^p(\calA) \ar@{->}[r]^{\frac{dz_1}{z_1}\ \ \ \ }\ar@{->}[d]_{\bfres_H^p}  & \frac{dz_1}{z_1}\wedge\Omega^p(\calA) \ar@{->}[dl]^{\mbox{\scriptsize Residue map}}\\
M^p &.  \\
}
$$
First we consider the Hilbert series of the submodules 
$(dz_1/z_1)\wedge\Omega^{p-1}(\calA)\subset \Omega^p(\calA)$. 
The following short exact sequence is obtained from Prop \ref{prop:exact}, 
$$
\begin{CD}
0 @>>> \frac{dz_1}{z_1}\wedge\Omega^{p-1}(\calA) 
@>>> \Omega^p(\calA)@>\frac{dz_1}{z_1}\wedge >> 
\frac{dz_1}{z_1}\wedge\Omega^{p}(\calA) @>>>0. 
\end{CD}
$$
We have a formula on Hilbert series ($p=0,1, \cdots, \ell$)
$$
P\left(\Omega^p(\calA), x\right)=
P\left( \frac{dz_1}{z_1}\wedge\Omega^{p-1}(\calA), x\right)+
x\cdot P\left( \frac{dz_1}{z_1}\wedge\Omega^{p}(\calA), x\right). 
$$
Summing up with multiply $y^p$, we have 
\begin{equation}
\label{eq:hilb}
\Phi(\calA ; x, y)
=\left(1+\frac{x}{y}\right)\times
\sum_{p=1}^\ell P\left( \frac{dz_1}{z_1}\wedge\Omega^{p-1}(\calA), x\right)y^p
\end{equation}
Next we consider the residue map. For general graded $S$-module $M$, 
suppose 
$s\in S$(with $\deg =1$) is an $M$-regular element, i.e. 
$s\cdot :M\stackrel{s}{\longrightarrow}M$ is injective, 
then from the exact sequence $0\rightarrow sM\rightarrow M
\rightarrow M/sM\rightarrow 0$, we can compute 
the Hilbert series of $M/sM$ as 
\begin{equation}
\label{eq:m/zm}
P(M/sM, x)=(1-x)P(M, x). 
\end{equation}
Using (\ref{eq:m/zm}), 
we can relate the Hilbert series of the module 
$M^p\subset\Omega^p\left(\calA^{H_1}, \bfk\right)$ with 
that of $\Omega^p(\calA)$. 
\begin{equation}
\label{eq:formulares}
P\left(M^p, x\right)=
x(1-x)P\left(\frac{dz_1}{z_1}\wedge\Omega^{p}(\calA), x\right). 
\end{equation}
Combining (\ref{eq:hilb}) and (\ref{eq:formulares}), 
we have the desired result. \owari

\section{Free arrangements in $\bbK^3$}
\label{sec:3arr}

In this section, let $\calA$ be a central arrangement in $V=\bbK^3$. 

Given a hyperplane $H\in\calA$, we have considered the natural 
multiarrangement $(\calA^{H}, \bfk)$ on $H$. 
Since $\Omega^1(\calA^H, \bfk)$ is reflexive $\bbK[H]$-module
and $\dim\bbK[H]=2$, 
$\Omega^1(\calA^H, \bfk)$ 
is a free $\bbK[H]$-module, we let $(d'_2, d'_3)$ denote the multiexponents. 
Since the sum of multiplicities of the multiarrangement  $(\calA^H, \bfk)$ 
is $\sharp(\calA)-1$, 
from Saito's criterion \ref{thm:saito's}, we have 
\begin{equation}
\label{eq:sum}
d'_2+d'_3=\sharp(\calA)-1. 
\end{equation}
We call $\chi_0(\calA, t):=\chi(\calA, t)/(t-1)$ 
the {\bf reduced characteristic polynomial}. 
If $\calA$ is a free arrangement, 
it follows from results of Ziegler (Theorem\ref{thm:zie}) that the 
exponents of $\calA$ are $\exp(\calA)=(1, d'_2, d'_3)$, and by 
Terao (Theorem \ref{thm:fact}) that 
\begin{equation}
\label{eq:coincide}
\chi_0(\calA, t)=(t-d'_2)(t-d'_3). 
\end{equation}
Conversely, if the multiexponent $\exp(\calA, \bfk)=(d_2', d_3')$ are not 
roots of the characteristic polynomial $\chi(\calA, t)$, $\calA$ cannot 
be free even if $\chi(\calA, t)$ factors into 
linear terms over $\bbZ$. 

\begin{Example}
\label{ex:stanley}{\normalfont (Stanley's example)}
Let $\calA$ be the cone of the affine $2$-arrangement 
in Fig.\ref{fig:stanley}(left). Then the characteristic 
polynomial factors as $\chi(\calA, t)=(t-1)(t-3)^2$. However, 
the restriction of $\calA$ to the hyperplane at infinity, 
Fig.\ref{fig:stanley}(right), 
has exponents $(1, 5)$, hence $\calA$ is not free. 
\end{Example}

\setlength{\unitlength}{.4pt}
\begin{figure}[htbp]
  \begin{center}
    \leavevmode
    \begin{picture}(200,200)(0,0)
      \put(250,0){\line(1,2){100}}
      \put(350,0){\line(-1,2){100}}
      \put(189,100){\line(1,0){223}}
      \put(300,0){\line(0,1){200}}
      \put(200,0){\line(1,1){200}}
      \put(200,200){\line(1,-1){200}}
      \put(300,100){\circle*{6}}

      \put(-150,0){\line(1,2){100}}
      \put(-50,0){\line(-1,2){100}}
      \put(-211,150){\line(1,0){223}}
      \put(-100,0){\line(0,1){200}}
      \put(-200,0){\line(1,1){200}}
      \put(-200,200){\line(1,-1){200}}
      \put(-100,100){\circle*{6}}
    \end{picture}
    \caption{Stanley's example and its restriction to the hyperplane at infinity.}
    \label{fig:stanley}
  \end{center}
\end{figure}
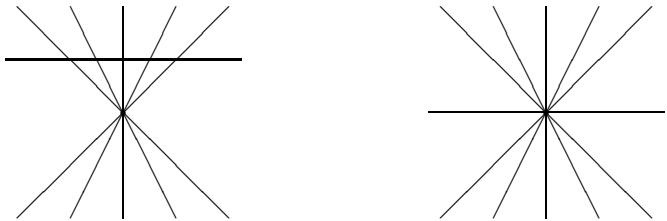

Because of relation (\ref{eq:sum}), the equation (\ref{eq:coincide}) 
is equivalent to 
$\chi_0(\calA, 0)=d'_2\cdot d'_3$. 
The main result 
of this section is that these relations characterize the 
freeness of $\calA$. 

\begin{Thm}
\label{thm:3}
Let $\calA$ be an arrangement in $\bbK^3$, 
$\chi_0(\calA, t):=(1-t)^{-1}\chi(\calA, t)$ be the reduced 
characteristic polynomial, $(d'_2, d'_3)$ be the multiexponents 
of restricted multiarrangement $(\calA^H, \bfk)$ and $M^1$ be the 
image of restriction map 
$\bfres_H^1:\Omega^1(\calA)\longrightarrow\Omega^1(\calA^H,\bfk)$. 
Then the codimension of $M^1$ in $\Omega^1(\calA^H,\bfk)$ is 
finite and is given by 
$$
\chi_0(\calA, 0)-d'_2\cdot d'_3. 
$$
\end{Thm}
In particular, from Corollary \ref{cor:criterion}, 
\begin{Cor}
\label{cor:char3}
$\calA$ is free if and only if 
$
\chi(\calA, t)=(t-1)(t- d'_2)(t-d'_3), 
$
where $(d_2', d_3')$ are the multiexponents of the 
restricted multiarrangement. 
\end{Cor}

\noindent
{\bf Proof of Theorem.} 
From the assumption, $\Omega^1(\calA^H,\bfk)$ has 
a homogeneous basis with degrees $(-d_2', -d_3')$. 
Hence it follows from the simple computation that 
\begin{equation}
\label{eq:simple}
\lim_{x\rightarrow 1}\sum_{p=0}^{2}
P(\Omega^p(\calA^H,\bfk), x)(t(1-x)-1)^p=(t-d_2')(t-d_3'). 
\end{equation}
Recall the characteristic polynomial can be calculated from 
$$
\Phi(\calA: x, y)=\sum_{p=0}^{3}P\left(\Omega^p(\calA), x\right)y^p 
$$
by Theorem \ref{thm:st}. We compare the Hilbert series above to 
that of $\Omega^p(\calA^H,\bfk)$ ($p=0,1,2$). From 
Theorem \ref{thm:hilbformula}, $\Phi(\calA; x, y)$ can be expressed 
by the Hilbert series of $M^p$ 
$$
\Phi(\calA; x,y)=
\frac{x+y}{x(1-x)}\sum_{p=0}^{2}P\left(M^p, x\right)
y^p. 
$$
Hence, 
\begin{eqnarray*}
\chi_0(\calA, t) &=& \frac{1}{t-1}\lim_{x\rightarrow1}
\frac{x+t(1-x)-1}{x(1-x)}\sum_{p=0}^{2}P\left(M^p, x\right)
(t(1-x)-1)^p \\
&=& \lim_{x\rightarrow1}
\sum_{p=0}^{2}P\left(M^p, x\right)(t(1-x)-1)^p . 
\end{eqnarray*}
We note that the maps $\bfres_H^0, \bfres_H^2$ are 
naturally surjective. 
Hence, (recall (\ref{eq:sum})) 
\begin{eqnarray*}
&&P(M^0, x)=P(\Omega^0(\calA^H,\bfk), x)=\frac{1}{(1-x)^2}\\
&&P(M^2, x)=P(\Omega^2(\calA^H,\bfk), x)=
\frac{x^{1-\sharp(\calA)}}{(1-x)^2}, 
\end{eqnarray*}
and we have from the definition, 
\begin{equation}
\dim_\bbK\left(\Omega^1(\calA^H,\bfk)/M^1\right)=
\lim_{x\rightarrow 1}\left[P(\Omega^1(\calA^H,\bfk), x)
-P(M^1, x) \right]
\end{equation}
Using above relations, 
\begin{eqnarray*}
\chi_0(\calA, 0)-d'_1\cdot d'_2 
&=& \chi_0(\calA, t)-(t-d'_1)(t-d'_2)\\
&=& \lim_{x\rightarrow 1}(t(1-x)-1)\left[ P(M^1, x)-
P(\Omega^1(\calA^H,\bfk), x)\right]\\
&=& t\cdot \lim_{x\rightarrow 1}(1-x)\left[ P(M^1, x)-
P(\Omega^1(\calA^H,\bfk), x)\right]\\
&&
\ \ \ \ \ -\lim_{x\rightarrow 1}\left[ P(M^1, x)-
P(\Omega^1(\calA^H,\bfk), x)\right]
\end{eqnarray*}
Since left hand side is a constant, we have 
$$
\dim_\bbK\left(\Omega^1(\calA^H,\bfk)/M^1\right)=
\chi_0(\calA, 0)-d'_1\cdot d'_2. 
$$
\owari

\section{Free arrangements in $\bbK^{\ell +1}\ (\ell\geq 3)$}
\label{sec:4}

Let $\calA$ be an essential arrangement 
in $V=\bbK^{\ell +1}\ (\ell\geq 3)$. Fix a hyperplane 
$H_0$ and coordinate system $(z_0, z_1, \cdots, z_\ell)$ such that 
$H_0=\{ z_0 =0\}$. Denote by $S:=\bbK[z_0, z_1, \cdots, z_\ell]$ as usual. 

Define $S$-submodules  $D_0(\calA)$ and $\Omega_0^1(\calA)$ 
of $D(\calA)$ and $\Omega^1(\calA)$, respectively, by 
\begin{eqnarray*}
&&D_0(\calA):=\{\delta \in D(\calA)\ |\ \delta z_0 =0\}\\
&&\Omega_0^1(\calA):=\{ \omega \in \Omega^1(\calA)\ |\ <E, \omega>=0\}, 
\end{eqnarray*}
where $E$ is the Euler vector field and $<,>$ is the inner product. 
We have the following splitting as $S$-modules: 
\begin{eqnarray*}
&&D(\calA)=S\cdot E\oplus D_0(\calA)\\
&&\Omega^1(\calA)=S\cdot\frac{dz_0}{z_0}\oplus \Omega_0^1(\calA). 
\end{eqnarray*}
The duality between $D(\calA)$ and $\Omega^1(\calA)$ implies 
that the modules $D_0(\calA)$ and $\Omega_0^1(\calA)$ are dual
each other, hence they are reflexive. Contrary to previous sections, 
we consider the module of vector fields $D(\calA)$(or $D_0(\calA)$) 
in this section. 
Ziegler's restriction map can be formulated as follows. 
\begin{equation}
\label{eq:restmap}
\begin{array}{ccc}
D_0(\calA)& \longrightarrow & D(\calA^{H_1}, \bfk)\\
&&\\
\delta&\longmapsto &\delta|_{z_0=0}. 
\end{array}
\end{equation}
So we have an exact sequence, 
\begin{equation}
\label{eq:modexact}
0\longrightarrow D_0(\calA)
\stackrel{z_0\cdot}{\longrightarrow}
D_0(\calA)
\longrightarrow
D(\calA^{H_0}, \bfk). 
\end{equation}
Theorem \ref{thm:zie} and Corollary \ref{cor:criterion} can be 
also stated as follows: 
\begin{Thm}
\label{thm:zie2}
$\calA$ is free with $\exp(\calA)=(e_0(=1), e_1, \cdots, e_\ell)$ 
if and only if $D(\calA^{H_0}, \bfk)$ is free with exponents 
$(e_1, \cdots, e_\ell)$ and the restriction map (\ref{eq:restmap}) 
is surjective for some $H_0\in\calA$. 
\end{Thm}

From now on we consider 
reflexive $\calO_{\bbP^\ell}$-module 
$\til{D_0(\calA)}$ 
on $\bbP^\ell$ rather than the graded $S$-module 
$D_0(\calA)$ \cite{har-alg}. The local structure of 
the coherent sheaf 
$\til{D_0(\calA)}$ can be 
described by the local structrure of arrangement $\calA$. 
More precisely, if we define the localization of $\calA$ at $x$ by 
$$
\calA_x:=\{ H\in\calA\ |\ H\ni x\ \}
$$
for $x\in V$, the stalk $\til{D(\calA)}_{\bar{x}}$ 
at $\bar{x}\in\bbP^\ell$ is isomorphic to that of 
$\til{D(\calA_x)}$. 
\begin{equation}
\label{eq:stalk}
\til{D(\calA)}_{\bar{x}}\cong \til{D(\calA_x)}_{\bar{x}}. 
\end{equation}
In particular, $\til{D(\calA)}$ is a locally free sheaf on 
$\bbP^\ell$ if and only if $\calA$ is locally free, i.e. 
$\calA_X$ is free for all $X\in L_\calA\backslash\{0\}$. 
(For details see Musta\c t\v a and Schenck \cite[Thm 2.3]{mus-sch})

From the discussion above, we have an exact sequence 
$$
0\longrightarrow \til{D_0(\calA)}(d-1)
\longrightarrow \til{D_0(\calA)}(d)
\longrightarrow \til{D(\calA^{H_0}, \bfk)}(d)
$$
($d\in\bbZ$). But the last homomorphism is not necessarily surjective. 
For the sake of surjectivity we consider the following condition. 
\begin{Assum}
\label{assum}
Arrangement $\calA$ is locally free along $H_0$, i.e. 
$\calA_x$ is free for all $x\in H_0\backslash\{0\}$. 
\end{Assum}
Note that locally free arrangements satisfy this condition. 
\begin{Prop}
\label{prop:sheafexact}
If Condition \ref{assum} holds then the restriction map induces 
an exact sequence,  
\begin{equation}
0
\longrightarrow \til{D_0(\calA)}(d-1)
\longrightarrow \til{D_0(\calA)}(d)
\longrightarrow \til{D(\calA^{H_0}, \bfk)}(d)
\longrightarrow 0, 
\end{equation}
and we have 
$\til{D_0(\calA)}(d)|_{\bbP(H_0)}
=\til{D(\calA^{H_0}, \bfk)}(d)$, where 
$\bbP(H_0)$ is the projective hyperplane defined by $H_0\subset V$. 
\end{Prop}
{\bf Proof.} 
What we have to show is the surjectivity. 
Since $\calA_x$ is free for all $x\in H_0\backslash\{0\}$, 
the induced restriction map 
$$
D_0(\calA_x)\longrightarrow D(\calA_x^{H_0}, \bfk|_{\calA_x^{H_0}})
$$
is surjective by Theorem \ref{thm:zie2}. 
From (\ref{eq:stalk}), this shows the surjectivity 
of the homomorphism as sheaves. 
\owari

In the context of coherent sheaves on projective space, 
an arrangement $\calA$ is free with 
$\exp(\calA)=(e_0(=1), e_1, \cdots, e_\ell)$ 
if and only if the coherent sheaf 
$\til{D_0(\calA)}$ splits into a direct sum of line bundles as 
$$
\til{D_0(\calA)}=\calO_{\bbP^\ell}(-e_1)\oplus
\cdots\oplus\calO_{\bbP^\ell}(-e_\ell). 
$$
In the theory of vector bundles on projective space, the following 
remarkable theorem is known. 
\begin{Thm}
\label{thm:oss}
{\normalfont (\cite[Theorem 2.3.2.]{oss})}\\
A holomorphic vector bundle $\calE$ on $\bbP^\ell$ splits into a 
direct sum of line bundles precisely when its restriction to 
some plane $\bbP^2\subset\bbP^\ell$ splits. 
\end{Thm}
In particular, in the case $\ell\geq 3$, a vector bundle $\calE$ 
on $\bbP^\ell$ splits 
if its restriction to some hyperplane $\bbP^{\ell -1}$ splits. 
Here we consider the following condition. 
\begin{Assum}
\label{ass:mf}
The restricted multiarrangement  $(\calA^{H_0}, \bfk_{\calA}^{H_0})$ is 
free. 
\end{Assum}
From the discussion above, we conclude the following result. 
\begin{Cor}
Let $\calA$ be an arrangement in $\bbK^{\ell+1}\ (\ell\geq 3)$. 
$\calA$ is free if and only if $\calA$ is locally free and the 
restricted multiarrangement is free (=Condition \ref{ass:mf}) 
for some hyperplane of $\calA$. 
\end{Cor}
The aim of this section is to characterize freeness 
by Condition \ref{assum}. This condition is 
weaker than local freeness. 
If $\calA$ is not a locally free arrengement, 
$\til{D_0(\calA)}$ is not a vector bundle. 
In this case the proof of Theorem \ref{thm:oss} does not work 
for $\til{D_0(\calA)}$. 
We recall some results on reflexive sheaves 
over projective spaces. For details see \cite{har-ref}. 

A coherent sheaf $\calF$ on a scheme $X$ is said to be reflexive 
if the natural map $\calF\rightarrow\calF^{\lor\lor}$ of $\calF$ to 
its double dual is an isomorphism. If $X$ is an regular scheme, i.e. 
all its local rings are regular local rings, then a reflexive sheaf 
$\calF$ on $X$ is locally free except along a closed subscheme $Y$ 
of codimension $\geq 3$ (\cite[Corollary 1.4.]{har-ref}). 

Now we prove 
\begin{Lemma}
\label{lem:refl}
Let $\calF$ be a reflexive sheaf on $\bbP^n\ (n\geq 3)$ and 
assume that $\calF$ is locally free except at a finite number of 
points $\{P_i\}$. Then there exist a surjection 
$$
\rmH^{n-1}(\calF^\lor\otimes \omega)\longrightarrow 
\rmExt^{n-1}(\calF, \omega)\longrightarrow 0, 
$$
where $\omega=\calO(-n-1)$ is the dualizing sheaf on $\bbP^n$. 
\end{Lemma}

\noindent
{\bf Proof.} (See also \cite[Theorem 2.4.]{har-ref}) 
Consider the spectral sequence of local and global Ext functors: 
$$
E_2^{p,q}
=\rmH^p(\calExt^q(\calF, \omega))\Longrightarrow
E^{p+q}=\rmExt^{p+q}(\calF, \omega). 
$$
From the assumption, the non-zero $E_2^{p,q}$ terms appear 
either $p=0$ or $q=0$. Indeed, since $\calF$ is 
locally free except for $\{P_i\}$, supports of 
$\calExt^q(\calF, \omega)\ (q\geq 1)$ are contained in 
finite set $\{P_i\}$. 
It follows that $E_2^{p,q}=0$ for $p,\ q\geq 1$. 
Furthermore, at these points $\calF$ has depth $\geq 2$
(\cite[Prop. 1.3.]{har-ref}), hence has 
homological dimension $\leq n-2$. Thus 
$\calExt^q(\calF, \omega)=0$ for $q\geq n-1$. Now the 
result is straightfoward from the definition of 
spectral sequence. 
\owari

Since $\rmExt^{n-1}(\calF(d), \omega)$ is Serre dual 
to $\rmH^1(\bbP^n, \calF(d))$, we have the inequality 
$$
\dim \rmH^1(\calF(d))\leq \dim\rmH^{n-1}(\calF(d)^\lor\otimes\omega), 
$$
by the theorem. 
The right-hand side vanishes for $d\ll 0$. 
\begin{Cor}\label{cor:vanish}\ 
$\rmH^1(\bbP^n, \calF(d))=0$ for $d\ll 0$. 
\end{Cor}

We now prove the main theorem of this section. 
\begin{Thm}
\label{thm:FREE}
Let $\calA$ be a arrangement in $\bbK^{\ell+1}\ (\ell\geq3)$ and 
fix a hyperplane $H\in\calA$. 
Then $\calA$ is free if (and only if) $\calA$ is locally free 
along $H$ 
(=Condition \ref{assum}) and the restricted multiarrangement 
$(\calA^{H}, \bfk_\calA^{H})$ 
is free (=Condition \ref{ass:mf}). 
\end{Thm}

\noindent
{\bf Proof.} We first note that $\til{D_0(\calA)}$ is 
locally free except for finite points. Indeed, if there exists 
an $X\in L_\calA$ with $\dim X\geq 2$ 
such that $\calA_X$ is not free, 
$X$ must intersect with the hyperplane $H$. 
Then $X\cap H$ is a set at which $\til{D_0(\calA)}$ is not locally free, 
which contradicts 
the assumption that $\calA$ is locally free along $H$. 

From the vanishing of intermediate 
cohomology groups of line bundles over projective space, 
we have 
$$
\rmH^i(\bbP(H), \til{D_0(\calA)}|_{\bbP(H)})=0, 
\mbox{ for }1\leq i\leq \ell -2. 
$$
Hence the next exact sequence is obtained
($\calF:=\til{D_0(\calA)}$). 
$$
\begin{array}{ccccccc}
0
&\longrightarrow&\rmH^0(\calF(d-1))
&\longrightarrow&\rmH^0(\calF(d))
&\longrightarrow&\rmH^0(\calF(d)|_{\bbP(H)})\\
&&&&&&\\
&\longrightarrow&\rmH^1(\calF(d-1))
&\longrightarrow&\rmH^1(\calF(d))
&\longrightarrow&0. 
\end{array}
$$
The surjection in the second row and Corollary \ref{cor:vanish} indicate 
that $\rmH^1(\calF(d))=0$ for any $d\in\bbZ$. Thus 
$\rmH^0(\calF(d))\rightarrow\rmH^0(\calF(d)|_{\bbP(H)})$ is surjective 
for any $d\in\bbZ$. This implies that 
$$
D_0(\calA)\rightarrow D(\calA^H, \bfk_\calA^H)\rightarrow 0
$$
is surjective. From (\ref{cor:criterion}) we conclude that $\calA$ is free. 
\owari

\section{Application}

We use Theorem \ref{thm:FREE} to show that cones over a certain 
truncated affine Weyl arrangements are free. 

Let $V=\bbR^\ell$ be an $\ell$-dimensional Euclidean space 
with a coodinate system $(x_1, \cdots, x_\ell)$. 
Let $\Phi$ be a crystallographic irreducible root system in $V$ 
with exponents $(e_1, e_2, \cdots, e_\ell)$ and Coxeter number $h$. 
We also fix 
a positive root system $\Phi^+\subset \Phi$. 
For each positive root $\alpha\in\Phi^+$ and integer $k\in\bbZ$, 
define an affine hyperplane $H_{\alpha, k}$ by 
$$
H_{\alpha,k}:=\{v\in V\ |\ \alpha(v) =k\}.
$$ 
We have a hyperplane arrangement 
$$
\calA(\Phi^+):=\{ H_{\alpha,0}\ |\ \alpha \in\Phi^+\ \}
$$
in $V$, called the Weyl arrangement associated to $\Phi$. 
The Weyl arrangement is free, more generally, 
\begin{Thm}
{\normalfont \cite{sai-linear}} 
Let $\calA$ be a Coxeter arrangement, i.e. the set of all reflecting 
hyperplanes of a finite Coxeter group of exponents $(e_1, \cdots, e_\ell)$. 
Then $\calA$ is a free arrangement with $\exp(\calA)=(e_1, \cdots, e_\ell)$. 
\end{Thm}
The basis of $D(\calA)$ can be constructed explicitly 
in terms of the invariant theory of Coxeter groups. 
We next define a family of arrangements in $\bbR\times V$ (with coordinate 
system $(x_0, x_1, \cdots, x_\ell)$) associated with an affine 
Weyl arrangement. These kinds of arrangements were first studied by 
Shi \cite{shi-lmn, shi-sign}. 
\begin{Def}
For integers $p, q\in\bbZ$ with $p\leq q$, denote by $[p, q]$ the set 
$\{p, p+1, \cdots, q\}$ of integers from $p$ to $q$. We define
an affine arrangement in $V$ as follows: 
$$
\calA(\Phi^+)^{[p,q]}:=\{ H_{\alpha, k}\ |\ \alpha\in\Phi^+,\ k\in [p, q]\ \}, 
$$
and its cone in $\bbR\times V$ 
$$
\cone\calA(\Phi^+)^{[p, q]}:=\{ \alpha -kx_0\ |\ \alpha\in\Phi^+,\ k\in [p, q]\ \}
\cup\left\{H_\infty:=\{x_0=0\}\right\}. 
$$
\end{Def}

Edelman and Reiner \cite{ede-rei} posed a conjecture on the freeness 
for such kind of arrangements. 
\begin{Conjecture}{\normalfont \cite{ede-rei}}
\label{conj:er}
\begin{itemize}
\item[{(1)}] The cone $\cone\calA(\Phi^+)^{[-m, m]}\ (m\geq 0)$ 
of an extended Catalan arrangement 
$\calA(\Phi^+)^{[-m, m]}$, $(m\geq 0)$ 
is free with exponents $(1, e_1+mh, e_2+mh, \cdots, e_\ell+mh)$. 
\item[{(2)}] The cone $\cone\calA(\Phi^+)^{[1-m, m]}\ (m\geq 1)$ 
of an extended Shi arrangement 
$\calA(\Phi^+)^{[1-m, m]}$, $(m\geq 1)$ 
is free with exponents $(1, mh, mh, \cdots, mh)$. 
\end{itemize}
\end{Conjecture}
By the general theory of free arrangements, we can deduce 
some conclusions from the conjecture. First restricting 
to the hyperplane at infinity $H_\infty=\{x_0=0\}$, 
from Theorem \ref{thm:zie}, we have. 
\begin{Conclusion}
\label{con:multi}
{\normalfont (Multi-freeness)}
For $\calA=\calA(\Phi^+)$ and $n\geq 0$, 
$$
D(\calA, n):=\{\delta\in\rmDer_V\ |\ \delta\alpha_H\in(\alpha_H^n),\ 
\forall H\in\calA\}
$$
is free with multiexponents 
$$
\left\{
\begin{array}{ll}
(e_1+mh, \cdots, e_\ell+mh),&\mbox{if}\ n=2m+1, \\
(mh, \cdots, mh),& \mbox{if}\  n=2m. 
\end{array}\right.
$$
\end{Conclusion}
Second, from the factorization 
theorem \ref{thm:fact}, we have another conclusion. 
\begin{Conclusion}
\label{con:fact}
{\normalfont (Factorization)}
 The characteristic polynomial of these arrangements 
are given by 
\begin{eqnarray*}
&&\chi(\calA(\Phi^+)^{[-m, m]}, t)=\prod_{i=1}^{\ell}(t-e_i -mh)\\
&&\chi(\calA(\Phi^+)^{[1-m, m]}, t)=(t-mh)^\ell. 
\end{eqnarray*}
\end{Conclusion}

In the case of root system of type $A$, 
Conjecture \ref{conj:er} (1) and (2) have been 
proved by Edelman and Reiner \cite{ede-rei} and Athanasiadis \cite{ath-a}, 
respectively. 

Without assuming conjecture, 
conclusion \ref{con:multi} was first studied by 
Solomon and Terao \cite{double}. 
Terao \cite{ter-multi} has completed the proof of 
(\ref{con:multi}) for Coxeter arrangements. 
A generalized version is also proved \cite{yos}. 
\begin{Thm}
\label{thm:multifree}
Let $\calA$ be a Coxeter arrangement with 
Coxeter number $h$. 
Suppose $\calA'\subset\calA$ is a free subarrangement with 
$\exp(\calA')=(e_1', \cdots, e_\ell')$. 
Let $\bfk$ be a multiplicity on $\calA$ defined by 
$$
\bfk(H)=
\left\{
\begin{array}{ll}
2m+1,&H\in\calA' \\
2m, &H\in\calA\backslash\calA'. 
\end{array}\right.
$$
Then the multiarrangement $(\calA, \bfk)$ is free with 
exponents 
$$
\exp(\calA, \bfk)=(e_1'+mh, \cdots, e_\ell'+mh)
$$
\end{Thm}

Conclusion \ref{con:fact} has been checked by computing 
characteristic polynomials 
when $\Phi$ is of classical type, that is 
of type $A$, $B$, $C$ or $D$, 
by Athanasiadis, see \cite{ath-def}. 
Recently Athanasiadis \cite{ath-cat} gives a case-free proof of 
the equation 
$$
\chi(\calA(\Phi^+)^{[-m, m]}, t)=\chi(\calA(\Phi^+), t-mh)
$$
which verifies \ref{con:fact} for extended Catalan arrangement. 
However, what we will need in our proof is a very weak version: 
(see also Prop.\ref{prop:gen} 
below) 
\begin{equation}
\label{eq:rk2}
\mbox{``Conclusion \ref{con:fact} is true for $A_2$, $B_2$ and $G_2$.''}
\end{equation}

Now let us prove Conjecture \ref{conj:er}. 
Our proof relies on the following elementary fact on root systems. 
\begin{Lemma}
Let $\Phi$ be an irreducible root system in $V$. Then for any point 
$x\in V\backslash\{0\}$, 
the localization of $\Phi$ at $x\in V$ 
$$
\Phi_x:=\{\alpha\in\Phi^+\ |\ \alpha (x)=0\ \}
$$
decomposes into a direct sum of root systems 
of lower ranks. 
\end{Lemma}
Let $\Phi_x=\Phi_1\oplus\cdots\oplus\Phi_k$ be an irreducible decompositon, 
a positive system $\Phi^+$ naturally determines 
positive systems on direct summands, by $\Phi_i^+:=\Phi^+\cap\Phi_i$. 

\begin{Thm}
\label{thm:erconj}
Conjecture \ref{conj:er} is true for any irreducible root system. 
\end{Thm}

\noindent
{\bf Proof.} We prove by induction on the rank of the root system $\Phi$. 
In the case of rank two, the cone $\cone\calA(\Phi^+)^{[p,q]}$ is a 
$3$-arrangement. Note that rank two root system is 
$A_2$, $B_2$ or $G_2$. 
(\ref{eq:rk2}) and 
Theorem \ref{thm:multifree} verifies the assumptions of 
Corollary \ref{cor:char3}. 
Thus the assertion is true for rank two root systems. 

Let $\Phi$ be an irreducible root system of higher rank. 
We apply Theorem \ref{thm:FREE} to prove freeness of 
$\cone\calA(\Phi^+)^{[p,q]}$, where $[p,q]$ is either $[1-m, m]$ or $[-m,m]$. 
The freeness of the restricted multiarrangement 
is verified by Theorem \ref{thm:multifree}.  
So what we have to check is 
that $\cone\calA(\Phi^+)^{[p,q]}$ 
is locally free along $H_\infty$. 

$H_\infty$ can be identified with 
$V\cong \{0\}\times V$. 
For given point $x\in V\backslash\{0\}$, put $\Phi_x=\Phi_1\oplus\cdots\oplus\Phi_k$, 
where $\Phi_i$ are irreducible root systems which have strictly 
lower rank than 
that of $\Phi$. 
Then it is easily seen that the localization of 
$\cone\calA(\Phi^+)^{[p,q]}$ at $x\in H_\infty\backslash\{0\}$ is 
$$
(\cone\calA(\Phi^+)^{[p,q]})_x=\cone\left(\calA(\Phi_1^+)^{[p,q]}
\oplus\cdots\oplus\calA(\Phi_k^+)^{[p,q]}\right). 
$$
From the inductive assumption, $\cone\calA(\Phi_i^+)^{[p,q]}$ 
is free for each $i=1, \cdots, k$. 
To finish the proof, we apply the next lemma. 
\owari 

\begin{Lemma}
Let $\calA_1,\cdots, \calA_k$ be affine arrangements such that 
each cone $\cone\calA_i$ is free ($i=1, \cdots, k$). Then 
$\cone(\calA_1\oplus\calA_2\oplus\cdots\oplus\calA_k)$ is also free. 
\end{Lemma}

As a possible generalization of Edelman-Reiner conjecture, 
we give a family of free arrangements interpolating 
between extended Shi $\cone\calA(\Phi^+)^{[1-m,m]}$
and Catalan $\cone\calA(\Phi^+)^{[-m,m]}$ arrangements. 
Recall that for positive roots $\alpha, \beta\in\Phi^+$, 
we denote $\alpha\leq\beta$ if $\beta -\alpha$ is 
a nonnegative linear combination of simple roots. 
Let $\Psi\subset\Phi^+$ be a subset of positive roots 
satisfying following conditions (1) and (2). 
\begin{itemize}
\item[(1)] $\Psi\subset\Phi^+$ is an order ideal, i.e. 
$\alpha\in\Psi$ and $\beta\leq\alpha\Longrightarrow\beta\in\Psi$. 
\item[(2)] $\calA(\Psi):=\{H_{\alpha, 0} |\ \alpha\in\Psi \}$ is a 
free arrangement (letting $\exp(\calA(\Psi))=(e_1', e_2', \cdots, e_\ell')$). 
\end{itemize}
\begin{Rem}
\normalfont
We do not know if (1) implies (2). 
When $\Phi$ is of type $A_{n-1}$, a subarrangement of $\calA(\Phi)$ 
corresponds to a graph with $n$ vertices, which is called a 
graphic arrangements. In this case, 
by using Stanley's characterization of 
freeness of graphic arrangement \cite{sta-sup} (see also \cite{between}), 
any subarrangement $\calA(\Psi)$ 
determined by an order ideal $\Psi\subset\Phi^+$ is free. 
\end{Rem}
For $m\in\bbZ_{\geq 0}$, let us take a union of 
the extended Shi $\calA(\Phi^+)^{[1-m,m]}$ and $\Psi$, 
$$
\calA(\Phi^+, \Psi, m):=\calA(\Phi^+)^{[1-m, m]}\cup
\{H_{\alpha, -m}\ |\ \alpha\in\Psi\}. 
$$
For our generalization, (\ref{eq:rk2}) should be generalized as 
follows, which can be proved by elementary computations. 
\begin{Prop}
\label{prop:gen}
Let $\Phi_0$ be a root system of rank two, $\Psi_0\subset\Phi_0^+$ an 
order ideal and $h$ the Coxeter number. 
Denote by $(e_1', e_2')=(1, \sharp(\Psi_0) -1)$ the 
exponents of $\calA(\Psi_0)$. Then the characteristic polynomial 
of $\calA(\Phi_0^+, \Psi_0, m)$ is 
$$
\chi(\calA(\Phi_0^+, \Psi_0, m), t)=(t-e_1'-mh)(t-e_2'-mh). 
$$
\end{Prop}
The proof of 
the next theorem is so similar to that of Theorem \ref{thm:erconj} 
that it will be omitted. 
\begin{Thm}
\label{thm:gener}
With notation as above, 
the cone $\cone\calA(\Phi^+, \Psi, m)$ 
is free with exponents
$$
\exp(\cone\calA(\Phi^+, \Psi, m))=(1, e_1'+mh, e_2'+mh, \cdots, e_\ell'+mh). 
$$
\end{Thm}
\begin{Cor}
$\chi\left(\calA(\Phi^+, \Psi, m), t\right)=\chi(\calA(\Psi), t-mh). $
\end{Cor}

\setlength{\unitlength}{.4pt}
\begin{figure}[htbp]
  \begin{center}
    \leavevmode
    \begin{picture}(200,200)(0,0)
      \put(-100,100){\circle*{10}}
      \put(-100,100){\vector(0,1){80}}
      \put(-100,100){\vector(2,1){68}}
      \put(-100,100){\vector(2,-1){68}}

      \put(-211,100){\line(1,0){223}}

\multiput(-100,100)(3,6){18}{\circle*{2}}
\multiput(-100,100)(-3,6){18}{\circle*{2}}
\multiput(-100,100)(3,-6){18}{\circle*{2}}
\multiput(-100,100)(-3,-6){18}{\circle*{2}}

\put(-230, 90){$\Psi$}


      \put(300,100){\circle*{10}}

      \put(250,0){\line(1,2){100}}
      \put(300,0){\line(1,2){100}}

      \put(350,0){\line(-1,2){100}}
      \put(400,0){\line(-1,2){100}}

      \put(189,150){\line(1,0){223}}
      \put(189,100){\line(1,0){223}}
      \put(189,50){\line(1,0){223}}
    \end{picture}
    \caption{Subarrangement $\Psi$ and $\calA(A_2, \Psi, 1)$}
    \label{fig:zu}
  \end{center}
\end{figure}
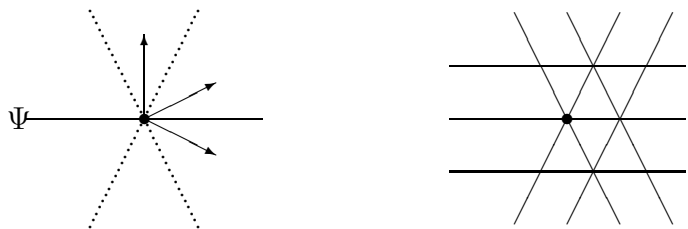

Note that this family contains both extended Shi and Catalan arrangements. 
In fact, taking $\Psi=\Phi^+$ we have the extended Catalan arrangement 
$\calA(\Phi^+, \Phi^+, m)=\calA(\Phi^+)^{[-m, m]}$ and 
taking $\Psi=\phi$(empty arrangement), we have the extended Shi 
arrangement 
$\calA(\Phi, \phi, m)=\calA(\Phi^+)^{[1-m, m]}$. 




\begin{thebibliography}{99}

\bibitem[Ath1]{ath-adv}
C. A. Athanasiadis, 
Characteristic polynomials of subspace arrangements and finite fields. 
{\it Adv. Math.} {\bf 122}  (1996), no. 2, 193--233. 

\bibitem[Ath2]{ath-a}
C. A. Athanasiadis, 
On free deformations of the braid arrangement. 
{\it European J. Combin.} {\bf 19} (1998), no. 1, 7--18.

\bibitem[Ath3]{ath-def}
C. A. Athanasiadis, 
Deformations of Coxeter hyperplane arrangements and their characteristic polynomials, 
in {\it Arrangements - Tokyo 1998}, pp.1-26, 
Advanced Studies in Pure Mathematics {\bf 27}, Kinokuniya, Tokyo, 2000. 

\bibitem[Ath4]{ath-cat}
C. A. Athanasiadis, 
Generalized Catalan numbers, Weyl groups and arrangements of hyperplanes. 
Preprint 

\bibitem[ER1]{between}
P. H. Edelman and V. Reiner, 
Free hyperplane arrangements between $A_{n-1}$ and $B_n$. 
{\it Math. Z. } {\bf 215} (1994), 347--365. 

\bibitem[ER2]{ede-rei}
P. H. Edelman and V. Reiner, 
Free arrangements and rhombic tilings. 
{\it Discrete Comput. Geom.} {\bf 15} (1996), no. 3, 307--340. 

\bibitem[Ha1]{har-alg}
R. Hartshorne, 
Algebraic geometry. Graduate Texts in Mathematics, No. 52. Springer-Verlag, New York-Heidelberg,  1977. 

\bibitem[Ha2]{har-ref}
R. Hartshorne, 
Stable reflexive sheaves. 
{\it Math. Ann.} {\bf 254} (1980), no. 2, 121--176. 

\bibitem[Ku]{kun}
J. P. S. Kung, 
A geometric condition for a hyperplane arrangement to be free. 
{\it Adv. Math.} {\bf 135} (1998), no. 2, 303--329.

\bibitem[MS]{mus-sch}
M. Musta\c t\v a and H. Schenck, 
The module of logarithmic $p$-forms of a locally free arrangement. 
{\it J. Algebra}  {\bf 241} (2001), no. 2, 699--719. 

\bibitem[OSS]{oss}
C. Okonek, M. Schneider and H. Spindler, 
Vector bundles on complex projective spaces. Progress in Mathematics,
  3. Birkh\"auser, Boston, Mass., 1980. vii+389 pp.

\bibitem[OT]{orl-ter}
P. Orlik and H. Terao, 
Arrangements of hyperplanes. 
Grundlehren der Mathematischen Wissenschaften, 300. 
Springer-Verlag, Berlin, 1992

\bibitem[Sa1]{sai-log}
K. Saito, 
Theory of logarithmic differential forms and logarithmic vector fields. 
{\it J. Fac. Sci. Univ. Tokyo Sect. IA Math.}{\bf 27} (1980), no. 2, 265--291

\bibitem[Sa2]{sai-linear}
K. Saito, 
On a linear structure of the quotient variety by a finite reflexion group. 
{\it Publ. Res. Inst. Math. Sci.} {\bf 29} (1993), no. 4, 535--579.

\bibitem[Sc1]{sch-rk2}
H. K. Schenck, 
A rank two vector bundle associated to a three arrangement, and its Chern polynomial. 
{\it Adv. Math.}  {\bf 149} (2000), no. 2, 214--229.

\bibitem[Sc2]{sch-recent}
H. K. Schenck, 
Elementary modifications and line configurations in $\bbP^2$.
 Comment. Math. Helv. {\bf 78} (2003), no. 3, 447--462.

\bibitem[Sh1]{shi-lmn}
J.-Y. Shi, 
The Kazhdan-Lusztig cells in certain affine Weyl groups. 
Lecture Notes in Math., {\bf 1179}, Springer Verlag, 1986 

\bibitem[Sh2]{shi-sign}
J.-Y. Shi, 
Sign types corresponding to an affine Weyl group. 
{\it J. London Math. Soc.} {\bf 35} (1987) 56--74 

\bibitem[ST1]{sol-ter}
L. Solomon and H. Terao, 
A formula for the characteristic polynomial of an arrangement. 
{\it Adv. in Math.} {\bf 64} (1987), no. 3,  305--325. 

\bibitem[ST2]{double}
L. Solomon and H. Terao, 
The double Coxeter arrangement. 
{\it Comment. Math. Helv.}{\bf 73} (1998), no. 2, 237--258.

\bibitem[Sta1]{sta-sup}
R. P. Stanley, 
Supersolvable lattices. 
{\it Algebra Universalis} {\bf 2} (1972), 197--217. 

\bibitem[Sta2]{sta-procnat}
R. P. Stanley, 
Hyperplane arrangements, interval orders, and trees. 
{\it Proc. Nat. Acad. Sci.}{\bf 93} (1996) 2620--2625. 

\bibitem[Sta3]{sta-rota}
R. P. Stanley, 
Hyperplane arrangements, parking functions, and tree inversions. 
Mathematical essays in honor of Gian-Carlo Rota 
(Cambridge, MA, 1996), 359--375, Progr. Math., 161, Birkh\"auser, 1998

\bibitem[Te1]{ter-arr}
H. Terao, 
Arrangements of hyperplanes and their freeness. I. 
{\it J. Fac. Sci. Univ. Tokyo Sect. IA Math.} {\bf 27} (1980), no. 2, 293--312.

\bibitem[Te2]{ter-free}
H. Terao, 
Generalized exponents of a free arrangement of hyperplanes and 
Shepherd-Todd-Brieskorn formula. 
{\it Invent. Math.} {\bf 63} (1981), no. 1, 159--179.

\bibitem[Te3]{ter-multi}
H. Terao, 
Multiderivations of Coxeter arrangements. 
{\it Invent. Math.} {\bf 148} (2002), no. 3, 659--674.

\bibitem[Yo]{yos}
M. Yoshinaga, 
The primitive derivation and freeness of multi-Coxeter arrangements. 
{\it Proc. Japan Acad. Ser. A Math. Sci.} {\bf 78} (2002), no. 7, 116--119.

\bibitem[Zi]{zie}
G. M. Ziegler, 
Multiarrangements of hyperplanes and their freeness, in 
{\it Singularities} (Iowa City,
IA, 1986), pp. 345--359, Contemp. Math., {\bf 90}, Amer. Math. Soc., Providence, RI, 1989. 


\end{thebibliography}
\end{document}